\pdfoutput=1
\documentclass[11pt,a4paper]{article}
\usepackage{fullpage}
\usepackage[USenglish]{babel}
\usepackage[T1]{fontenc}
\usepackage[utf8]{inputenc}
\usepackage{csquotes}
\usepackage{enumitem}
\usepackage{appendix}
\usepackage[giveninits,urldate=comp,sortcites,maxbibnames=10,url=false,backend=bibtex]{biblatex}
\usepackage[pdfencoding=auto,colorlinks=true,linkcolor=blue,bookmarksnumbered=true,pdfusetitle,hypertexnames=false]{hyperref}

\usepackage{amssymb}
\usepackage{amsmath}
\usepackage{mathtools}
\usepackage{amsthm}
\usepackage{dsfont}

\theoremstyle{plain}
\newtheorem{theo}{Theorem}[section]
\newtheorem{prop}[theo]{Proposition}

\newtheorem{lem}[theo]{Lemma}
\theoremstyle{definition}
\newtheorem{defin}[theo]{Definition}

\theoremstyle{remark}
\newtheorem{rem}[theo]{Remark}

\newcommand{\Ecal}{\mathcal E}

\newcommand{\Nds}{\mathds N}
\newcommand{\Qcal}{\mathcal Q}
\newcommand{\Rds}{\mathds R}
\newcommand{\Vbf}{\mathbf V}
\newcommand{\eps}{\varepsilon}
\newcommand{\wtd}{\widetilde}
\newcommand{\spt}{\mathrm{spt\,}}
\title{Lax--Oleinik formula on networks}
\author{Marco Pozza\thanks{Dipartimento di Matematica, Sapienza Universit\`a  di Roma, Italy. {\tt marco.pozza@uniroma1.it}} and Antonio Siconolfi\thanks{Dipartimento di Matematica, Sapienza Universit\`a  di Roma, Italy. {\tt siconolfi@mat.uniroma1.it}}}

\date{}

\begin{document}
\maketitle

\begin{abstract} We provide a Lax--Oleinik-type representation formula for solutions
of  time--dependent Hamilton--Jacobi equations, posed on a network with
a rather  general geometry, under standard assumptions on the Hamiltonians. It depends on a given initial datum at $t=0$ and a flux
limiter at the vertices, which both have to be assigned in order the problem to be uniquely solved.
Previous results in the same direction are solely in the frame of
junction, namely network with a single vertex. An important step  to get the result is
to   define a suitable action functional
and prove existence as well as  Lipschitz--continuity of minimizers between two fixed points
of the network in a given time, despite the fact that the integrand lacks convexity
  at the vertices.

\end{abstract}


\medskip


\section{Introduction}

The aim of the paper is to provide a Lax--Oleinik--type representation formula, obtained via
minimization of a suitable action functional,  for solutions
of  time--dependent Hamilton--Jacobi equations posed on a network with
a rather  general geometry.

This kind of formulas have actually a wider scope since they
 can enhance the  qualitative analysis of
Hamilton--Jacobi equations on networks, similarly to what happens for manifolds or
Euclidean spaces.  In this setting,  in fact, they essentially enter
 into play in a variety of theoretical constructions
such as    weak KAM theory, see \cite{fathi14},
 existence of regular subsolutions \cite{bernard07},  homogenization  problems \cite{souganidis99},
large time behavior of solutions \cite{fathi98,davinisiconolfi06,ishii08}, selection principles in the ergodic approximation \cite{davinifathiiturriaga16}.

Previous contributions on the same topic we deal with  can be found in \cite{ImbertMonneauZidani} where,
besides some restrictions on the Hamiltonians, the results are given  in the case of junctions,
namely networks with a single vertex. In \cite{IturriagaSanchezMorgado} the focus is instead
on the relationship between Lax--Oleinik formula and weak KAM theory and the connection
with time--dependent Hamilton--Jacobi problems are not taken into consideration.  Furthermore,
the Hamiltonians/Lagrangians are assumed of Tonelli type and required to satisfy a
quite stringent condition
at the vertices.

We consider a  connected network $\Gamma$ embedded in $\Rds^N$ with a finite number
of vertices, making up a set denoted by $\Vbf$,  linked by
regular simple curves $\gamma$ parametrized in $[0,1]$, called arcs of $\Gamma$.
A Hamiltonian on $\Gamma$ is defined  as a collection of one--dimensional Hamiltonians $H_\gamma:[0,1] \times \Rds \to \Rds$,
indexed by arcs, depending on state and momentum variable, with  the crucial feature
that Hamiltonians associated to arcs possessing different support, are totally unrelated.

We assume the $H_\gamma$'s  to  be continuous
in both arguments plus convex and superlinear in the momentum variable.
Namely, not more than  the usual conditions  to ensure the validity of
Lax--Oleinik formula for time--dependent  problems posed on  manifolds or  Euclidean spaces.

The  equations we are interested  in  are
\begin{equation}\label{intro1}
    U_t + H_\gamma(s,U')=0  \qquad\hbox{in $(0,1) \times (0,+\infty)$}
\end{equation}
on each arc $\gamma$,  and a solution on $\Gamma$ is   a continuous function
$u: \Gamma\times (0,+\infty) \longrightarrow \Rds$
such that $u(\gamma(s),t)$  solves in the viscosity sense \eqref{intro1} for each $\gamma$,
and satisfies suitable additional conditions on the discontinuity interfaces
\[\{(x,t), t \in [0,+\infty)\} \qquad\hbox{with $x \in \Vbf$.}\]
 It has been established in \cite{ImbertMonneau} in the case of junctions and
 in \cite{Siconolfi} for general networks that to get existence and uniqueness of solutions,
  equations \eqref{intro1} must be coupled not only with a continuous initial datum at $t=0$,
  but also  with a flux limiter, that is  a choice of appropriate  constants  $ c_x$ for $x$ varying in
   $\Vbf$.  We also report the contribution  of \cite{LionsSouganidis,Morfe}, where the time--dependent problem is studied in junctions,
    possibly multidimensional, with Kirchoff-type Neumann conditions at vertices.

In \cite{Siconolfi} flux limiters crucially appear in the conditions a solution must satisfy on the
   interfaces and, among other things, bond from above the time derivatives of any subsolution on it.
 Even if an initial datum is fixed,  solutions can change according to
 the choice of flux limiter,
  so that they must taken into account in representation formulas.

Due to the superlinearity of the $H_\gamma$'s, Lagrangians $L_\gamma$
can be given, for any arc,  via Fenchel transform. An overall  Lagrangian $L$, playing
the role of integrand  in the Lax--Oleinik formula,  is defined on
the whole tangent bundle  $T\Gamma$ of $\Gamma$ by
changing variables, gluing
together the $L_\gamma$'s and taking into account the chosen flux limiter $c_x$
via the formula
\begin{equation}\label{intro2}
 L(x,0)= c_x  \qquad\hbox{ for any $x \in \Vbf$.}
\end{equation}
 Setting $L=+\infty$ outside $T\Gamma$, we get a function which is lower semicontinuous in
$\Rds^N \times \Rds^N$, but with $L(x,\cdot)$ lacking in general convexity when $x$ is a vertex, see Section \ref{act1}.
 In fact the tangent space to $\Gamma$ at such points itself  is nonconvex being the union of
 one--dimensional vector spaces  corresponding to different intersecting arcs. This complicates
 to some extent the proof of the existence of curves minimizing the action functional
 between two fixed  points of $\Gamma$ in a given time, and of the
 Lipschitz  character of the minimizers as well, see   Theorems \ref{funda} and \ref{lippo}.

 In this respect condition \eqref{intro2} plays an important role since it
 allows ruling out Zeno phenomena in the minimization of
the action functional, namely the possibility that candidate minimizers wildly oscillate
around a discontinuity interface.

This is achieved taking into account the existence of regular solutions to the
 stationary equation $H_\gamma=-c_x$ in $(0,1)$, for arcs $\gamma$ ending at a given vertex  $x$, see  Proposition \ref{liscia},
in this way one can show that to oscillate  around the interfaces is more expensive
 than to stay steady at the vertex for a suitable time interval.
 This argument seems  new, see   Propositions \ref{postprecruxtris}  and \ref{precrux}.

In the end we get a formula similar to the one already known in the classical cases,
namely the solution $u(x,t)$ can be expressed as
\[u(x,t) = \inf \left \{\int_0^t L(\xi,\dot\xi) \, ds + u_0(\xi(0)) \right \}\]
with the infimum taken over the absolutely continuous curves $\xi$
  from $[0,t]$ to $\Gamma$ satisfying $\xi(t)=x$, where
$u_0$ stands for the initial datum and the flux limiter $c_x$ is taken into account through
the relation \eqref{intro2}.

The paper is organized as follows: in Section \ref{pre} we fix some notations and conventions.
Section \ref{net} provides  some basic facts about networks,   the definition of Hamiltonians/Lagrangians on networks and of flux limiters,  the main assumptions on the model are specified as well.  Section \ref{act} is devoted to the definition of the
overall Lagrangian $L$  and the corresponding action functional,
it is further introduced   a class of curves called admissible which is convenient to take
in consideration in the minimization of the action functional.
Existence of minimizing curves and their Lipschitz--continuity are proved in Section \ref{min}.

 Finally in  Section \ref{mai} it is proved  the continuity of the function given by Lax--Oleinik formula,  it is given the definition of solution for the time--dependent Hamilton--Jacobi equation on $\Gamma$, and written down    the main result stating  that Lax--Oleinik formula provides
 the unique solution of the  equation with  given initial datum and flux limiter.

Two appendices are about some elementary properties of  time--dependent equations on a
single arc, and  a variational result we use to prove the Lipschitz
 regularity of minimizers.

 \bigskip

\noindent {\bf Acknowledgement:}  The authors are members of INdAM--GNAMPA.

\noindent The authors would like to thank an anonymous referee for many remarks and suggestions on the manuscript which have allowed to improve the presentation and in particular  to essentially simplify  the proof of Theorem \ref{lippo} through the results presented in Appendix \ref{covidcovid}.

\bigskip

\section{Preliminaries} \label{pre}
We fix  a dimension $N$ and $\Rds^N$ as ambient space. We indicate by $\cdot$
the scalar product  in $\Rds^N$.
 Given two real numbers $a$ and $b$, we set
\[a \wedge b = \min \{a,b\}  \qquad\hbox{and} \qquad a \vee b = \max \{a,b\}.\]
  Given any subset $A$ contained in an Euclidean space, we denote by $\overline A$
its {\em closure}. Given a measurable  subset
$E \subset \Rds$, we denote by $|E|$ its {\em Lebesgue measure}. By curve contained   we mean throughout
the paper an {\em absolutely continuous} curve with support contained in $\Rds^N$ or $\Rds$.

\medspace

We set
 \[\Rds^+ = [0,+\infty),  \qquad  \Qcal = (0,1) \times (0,+\infty),\]
 we have $\partial \Qcal =[0,1]\times\{0\}\cup\{0,1\}\times\Rds^+$,
  we further set
  \[\partial^-\Qcal=[0,1]\times\{0\}\cup\{0\}\times\Rds^+.\]

  \medskip

For any $C^1$ function $\Phi: \Qcal \to \Rds$ and $(s_0,t_0) \in \Qcal$, we denote
by $\Phi'(s_0,t_0)$ the {\em space} derivative, with respect to $s$,  at
$(s_0,t_0)$, and by $\psi_t(s_0,t_0)$ or $\frac d{dt} \psi(s_0,t_0)$ the time derivative.

 Given a continuous function $U: \Qcal \to \Rds$, we call {\em supertangents}
 (resp.  {\em subtangents}) to $U$ at $(s_0,t_0) \in \Qcal$
 the viscosity test functions from above (resp. below).
 If needed, we take, without explicitly mentioning, $U$ and test function coinciding at
  $(s_0,t_0)$ and test function strictly greater (resp. less) than $U$ in a punctured
  neighborhood of $(s_0,t_0)$.

  We say that a  subtangent $\Phi$ to $U$ at $(1,t_0)$, $t_0 >0$, is
  {\em constrained to $\overline \Qcal$} if $(s_0,t_0)$ is  a minimizer of
  $u - \Phi$  in a neighborhood of $(1,t_0)$ intersected with $\overline Q$.

\bigskip

\section{Networks} \label{net}

\subsection{Basic definitions}  An {\it embedded network},  is a subset $
\Gamma \subset \Rds^N$ of the form
\[ \Gamma = \bigcup_{\gamma \in \Ecal} \, \gamma([0,1]) \subset \Rds^N,\]
where $\Ecal$ is a finite collection of regular ({\it i.e.}, $C^1$
with non-vanishing derivative) simple oriented curves, called {\it
arcs} of the network,  that we assume, without any loss of
generality, parametrized on $[0,1]$, note that we are also assuming existence of one--sided derivatives
at the endpoints $0$ and $1$. We stress out  that a regular change of parameter
does not affect our results.
 \\

On the support of any arc $\gamma$, we  also consider the
inverse parametrization   defined as
\[\widetilde \gamma(s)= \gamma( 1 -s) \qquad\hbox{for $s \in [0,1]$.}\]
We call $\widetilde \gamma$ the {\it inverse arc} of $\gamma$. We assume
\begin{equation}\label{netw}
    \gamma((0,1)) \cap \gamma'([0,1]) = \emptyset \qquad\hbox{whenever $\gamma \neq
\gamma'$, $\gamma \neq \widetilde\gamma'$}.
\end{equation}
In this case we say that $\gamma$, $\gamma'$ are {\em different arcs}.
\medskip
 We call  {\it vertices} the initial  and terminal points of the arcs,
 and denote  by  $\Vbf$ the
sets of all such vertices. Note that \eqref{netw} implies that
\[\gamma((0,1)) \cap \Vbf  = \emptyset \qquad\hbox{for any arc $\gamma$.}\]
 We assume that the network  is  {\em connected}, namely given two
vertices there is a finite concatenation of  arcs linking them. A {\em loop} is
an arc with initial and final point coinciding.
The unique  restriction we require  on the geometry of the network is
\begin{itemize}
  \item[{\bf (A1)}]  $\Ecal$ does not contain loops.
\end{itemize}
   See \cite{Sunada}
 for a comprehensive treatment on graphs and networks.\\

Given $x \in \Vbf$, we define
\[ \Gamma_x =\{ \gamma  \mid \gamma(1)=x\}.\]

\medskip

The network $\Gamma$ inherits a geodesic distance, denoted with $d_\Gamma$,
from the Euclidean metric of $\Rds^N$. It is clear that given $x$, $y$ in $\Gamma$ there is at least
a geodesic linking them.
The geodesic distance  is in addition equivalent to the Euclidean one.\\

 Given a continuous function $u: \Gamma  \times \Rds^+ \to \Rds$ and an arc $\gamma$,
we define $u \circ \gamma:[0,1] \times [0,+ \infty)
\to \Rds$ as
\begin{equation}\label{eq:uconcg}
u \circ \gamma(s,t)= u(\gamma(s),t) \qquad\hbox{for any $(s,t) \in \Qcal$.}
\end{equation}

\smallskip

The {\em tangent bundle} of $\Gamma$, $T \Gamma$ in symbols,
is made up by elements $(x,q) \in \Gamma \times \Rds^N$ with $q$ of the form
\[q = \lambda \, \dot \gamma(s)  \qquad\hbox{if $x=\gamma(s)$, $s \in [0,1]$,
with $\lambda \in \Rds$ },\]
   note that  $\dot \gamma(s)$ is univocally determined,
   up to a sign, if $x \in \Gamma \setminus \Vbf$ or in other words if $s \neq 0,1$.

\medskip

\subsection{Curves on \texorpdfstring{$\Gamma$}{Γ}}

We define

\begin{equation}\label{lippi1}
 m = \min\{|\dot\gamma(s)| \mid s \in [0,1],\, \gamma \in \Ecal \} >0.
\end{equation}

\smallskip

\begin{lem}\label{lippi} Given an arc $\gamma$, the function
\[\gamma^{-1}: \gamma[0,1] \to [0,1]\]
is Lipschitz continuous with respect to the distance $d_\Gamma$ on $\Gamma$.
\end{lem}
\proof It is enough to consider $x$, $y$ in $\gamma(0,1)$ with  $\gamma^{-1}(x) < \gamma^{-1}(y)$,
if the geodesic between $x$ and $y$ is given by the portion of $\gamma$ between
$\gamma^{-1}(x)$
and $\gamma^{-1}(y)$, we have
\[d_\Gamma(x,y) = \int_{\gamma^{-1}(x)}^{\gamma^{-1}(y)} |\dot\gamma| \, ds
\geq m \, (\gamma^{-1}(y)- \gamma^{-1}(x)).\]
If instead such a  geodesic is given by the concatenation of different arcs, then taking into account that the geodesic distance between vertices of $\Gamma$ is bounded by below by $m$, we get
\[d_\Gamma(x,y) \geq m \, (1-0) \geq m \, (\gamma^{-1}(y)- \gamma^{-1}(x)).\]
 We then derive in both cases
\[\gamma^{-1}(y)- \gamma^{-1}(x) \leq \frac 1m \, d_\Gamma(x,y).\]
\endproof

\smallskip

\begin{lem}\label{deri} For any given  arc $\gamma$  and   curve  $\xi:[a,b] \to \gamma([0,1])$,
 the function
\[\gamma^{-1} \circ \xi:[a,b] \to [0,1]\]
is absolutely continuous, and
\begin{equation}\label{deri1}
 \frac d{dt} \gamma^{-1} \circ \xi (t) = \frac{\dot\xi(t) \cdot \dot\gamma(\gamma^{-1} \circ \xi (t))}{|\dot\gamma(\gamma^{-1} \circ \xi (t))|^2}
   \quad\hbox{for a.e. $t$.}
\end{equation}
\end{lem}
\proof The function $\gamma^{-1} \circ \xi$ is absolutely continuous as composition of
 an absolutely continuous and a Lipschitz continuous function.
Let $t_0$ be a time where $\xi$ is differentiable, to ease notations we put
$s_0=\gamma^{-1} \circ \xi(t_0)$,
$\lambda_0= \frac{\dot \xi(t_0) \cdot \dot\gamma(s_0)}{|\dot \gamma(s_0)|^2}$ so that
\[\xi(t_0)= \gamma(s_0) \qquad\hbox{and}\qquad \dot \xi(t_0)= \lambda_0 \, \dot \gamma(s_0),\]
and formula \eqref{deri1} boils down to
\[
 \frac d{dt} \gamma^{-1} \circ \xi (t) = \lambda_0.
\]
We then have for $h \neq 0$
\begin{eqnarray}
 \frac{\gamma^{-1} \circ \xi(t_0+h) - \gamma^{-1} \circ \xi(t_0)}h &=& \frac{\gamma^{-1}(\xi(t_0) + h \, \dot \xi(t_0) + \mathrm{o}(h)) - \gamma^{-1}(\xi(t_0))}h  \nonumber\\
 &=&\frac{\gamma^{-1}(\gamma(s_0) + h \, \lambda_0 \, \dot\gamma(s_0) +\mathrm{o}(h))- s_0}h \label{deri12}\\
 &=& \frac{\gamma^{-1} (\gamma(s_0 + \lambda_0 \, h) + \mathrm{o}(h))- s_0}h,
\nonumber \end{eqnarray}
where $\mathrm o(\cdot)$ is the Landau symbol. Exploiting the Lipschitz continuity of $\gamma^{-1}$, we have
\[|\gamma^{-1} (\gamma(s_0 + \lambda_0 \, h) + \mathrm{o}(h))- \gamma^{-1} (\gamma(s_0 + \lambda_0 \, h))| \leq  \mathrm{o}(h)\]
or in other words
\[\gamma^{-1} (\gamma(s_0 + \lambda_0 \, h) + \mathrm{o}(h))=\gamma^{-1} (\gamma(s_0 + \lambda_0 \, h)) + \mathrm{o}(h).\]
Continuing the computation in \eqref{deri12} we therefore get
\begin{eqnarray*}
  \frac{\gamma^{-1} (\gamma(s_0 + \lambda_0 \, h) + \mathrm{o}(h))- s_0}h &=& \frac{\gamma^{-1} (\gamma(s_0 + \lambda_0 \, h)) + \mathrm{o}(h)- s_0}h \\
   &=& \frac{s_0+ \lambda_0 \,h -s_0 + \mathrm{o}(h) }h
   \\ &=& \frac {\lambda_0\,h + \mathrm{o}(h) }h.
\end{eqnarray*}
Sending $h$ to $0$, we obtain in the end
\[\lim_{h \to 0}\frac{\gamma^{-1} \circ \xi(t_0+h) - \gamma^{-1} \circ \xi(t_0)}h= \lambda_0 \]
which shows that $\gamma^{-1} \circ \xi$ is differentiable at $t_0$  and \eqref{deri1} holds true.
\endproof

\medskip

 \subsection{Hamiltonians and Lagrangians on the arcs} A Hamiltonian on  $\Gamma$ is a
 collection of Hamiltonians
$H_{\gamma}:[0,1] \times \Rds \to \Rds$, indexed by the arcs satisfying
\[
   H_{\widetilde\gamma}(s,\mu) = H_{\gamma}(1-s,-\mu) \qquad\hbox{for any arc $\gamma$} \\
\]
We emphasize that, apart the above compatibility condition, the Hamiltonians
$H_\gamma$ are {\em unrelated}.

\medskip

We  require  any  $H_\gamma$ to be:
\begin{itemize}
    \item[{\bf(H1)}]  continuous in both arguments;
    \item[{\bf(H2)}]  convex   in $\mu$;
    \item[{\bf(H3)}]  $\lim\limits_{|\mu|\to\infty}
\inf\limits_{s\in [0,1]}\dfrac{H_\gamma(s,\mu)}{|\mu|}= +\infty$
\; for any $\gamma\in\Ecal$.
   \end{itemize}

   \smallskip

Note that by the Corollary of Proposition 2.2.6 in \cite{Clarke}
the above assumptions imply that the $H_\gamma$'s are locally Lipschitz continuous  in $\mu$ uniformly with respect to $s \in [0,1]$.
     Namely, given $M >0$, there exists $C_M$   such that
\begin{equation}\label{newnew}
H_\gamma(s,\mu_1) - H_\gamma(s,q) \leq C_M \, |\mu_1-\mu_2|\qquad\hbox{for any $s \in [0,1]$, $\lambda_1$, $\lambda_2$ in $(-M,M)$}.
\end{equation}

  \smallskip

  Assumptions  {\bf (A1)}, {\bf(H1)}, {\bf(H2)}, {\bf(H3)} are in force,
  without further mentioning, throughout the paper.

  \smallskip

Thanks to the superlinearity condition {\bf (H3)}, we can define for any
 $\gamma\in\Ecal$,  the \emph{Lagrangian} corresponding  to  $H_\gamma$ as
\[L_\gamma(s,\lambda):=\max_{\lambda\in\Rds}(\lambda \, \mu -H_\gamma(s,\mu)).
\]
where the supremum is actually achieved thanks to {\bf (H3)}.
We have  for each $\lambda\in\Rds$ and $s \in [0,1]$,
\[
L_\gamma(s,\lambda)=L_{\widetilde\gamma}(1-s,-\lambda).
\]

From {\bf (H3)} we also derive that  the Lagrangians $L_\gamma$ are superlinear.
\medskip

We set
\[c_\gamma= - \max_s \, \min_\mu  H_\gamma(s,\mu) \qquad\hbox{for any arc $\gamma$.}\]

\begin{lem}\label{flusso} We have
\[c_\gamma = \min_s L_\gamma(s,0).\]
\end{lem}
\proof Given   $s \in [0,1]$, we have
\[L_\gamma(s,0)= \max_\mu -H_\gamma(s,\mu)= - \min_\mu H_\gamma(s,\mu)\]
and consequently
\[c_\gamma= - \max_s \, \min_\mu  H_\gamma(s,\mu) = \min_s [- \min_\mu  H_\gamma(s,\mu)] = \min_s L_\gamma(s,0).\]
\endproof

\smallskip

The following characterization  will play a relevant role in the sequel, see
 \cite{SiconolfiSorrentino}
for a comprehensive analysis of Eikonal equations on networks and graphs.

\smallskip

\begin{prop}\label{liscia} The stationary equation
\begin{equation}\label{liscia1}
  H_\gamma(s,U')=a  \qquad\hbox{ in $(0,1)$}
\end{equation}
admits a $C^1$ solution if and only if $a  \geq -c_\gamma$.
\end{prop}
\proof If $a$ satisfies the inequality in the statement, we define
\[\sigma^+_a(s)= \max\{\mu \mid H_\gamma(s,\mu)=a\},\]
which is apparently a continuous function as $s$ varies in $[0,1]$.
We further define
\[U(s)= \int_0^s \sigma_a^+(\tau) \, d\tau,\]
which is the sought $C^1$ solution to  \eqref{liscia1}.
The converse implication is immediate.
\endproof

\smallskip

Following \cite{ImbertMonneau}, we call {\em flux limiter} any
function  $x \mapsto c_x$ from $\Vbf$ to $\Rds$ satisfying
\[c_x  \leq \min_{\gamma\in \Gamma_x} c_\gamma = \min_\gamma L_\gamma(s,0) \qquad\hbox{for any $x \in \Vbf$.}\]

\smallskip

It is convenient for future use to introduce a modification of $L_\gamma$. For any arc $\gamma$ we define
 \[
   \bar L_\gamma(s,\lambda)=  \left \{\begin{array}{cc}
                               L_\gamma(s,\lambda) & \quad\hbox{for $s \neq 0,\,1$} \\
                                L_\gamma(s,\lambda) + (c_{\gamma(0)} \wedge c_{\gamma(1)}) - L_\gamma(s,0) & \quad\hbox{for $s = 0,\,1$}
                             \end{array} \right .
 \]

It is clear from the definition of $\bar L_\gamma$ and Lemma \ref{flusso}  that

\begin{lem}\label{postbarra} For any arc $\gamma$ the function $(s,\lambda) \mapsto
\bar L_\gamma(s,t)$ satisfies the following properties:
\begin{itemize}
  \item[{\bf (i)}] is lower semicontinuous in $(s,\lambda)$;
  \item[{\bf (ii)}] is convex in $\lambda$ for any fixed  $s \in [0,1]$;
  \item[{\bf (iii)}] $\bar L_\gamma(0,0)=\bar L_\gamma(1,0)= c_{\gamma(0)} \wedge c_{\gamma(1)}$;
  \item[{\bf (iv)}] is superlinear in $\lambda$ for any fixed $s$.
\end{itemize}
\end{lem}

\smallskip

The role of the flux limiters is highlighted  by the following result:

\begin{prop}\label{postprecruxtris} Given an arc $\gamma$ and a curve
 $\eta: [a,b] \to [0,1]$ with
   $\eta(a)= \eta(b)$,  one has
\[\int_a^b  L_\gamma(\eta,\dot\eta) \, d\tau \geq
\int_a^b \bar L_\gamma(\eta,\dot\eta) \, d\tau \geq (c_{\gamma(0)} \wedge
 c_{\gamma(1)}) \, (b-a).\]
\end{prop}
\proof
We set
\begin{eqnarray*}
   E &=& \{t \in [a,b] \mid \eta(t) \not \in \{0,\,1\}\}\\
   F &=& \{t \in [a,b] \mid \eta(t) \in \{0,\,1\}\}.
\end{eqnarray*}
If $E$ is empty then $\eta$ is constant, equal either to $0$ or to $1$, then
 the assertion is immediate in view of Lemma \ref{flusso}.
If  $E$ is instead nonempty,
it is open  in $[a,b]$, and it is then the union, up to a set of vanishing measure,
 of a countable family of intervals  $[a_k,b_k]$ with disjoint interiors. According to Proposition
 \ref{liscia}, there exists  a $C^1$ subsolution
$U(s)$ of $H_\gamma=- c$, where $c = c_{\gamma(0)} \wedge c_{\gamma(1)}$,
in $(0,1)$ so that we have
\[ \bar L_\gamma(\eta(t),\dot \eta(t))= L_\gamma(\eta(t),\dot \eta(t))
\geq \dot\eta(t) \, U'(\eta(t)) -H_\gamma(\eta(t), U'(\eta(t))) \]
for a.e. $t \in E$, and consequently
\begin{eqnarray*}
  \int_{a_k}^{b_k} \bar L_\gamma(\eta,\dot\eta) \, d\tau &\geq& \int_{a_k}^{b_k}
  U'(\eta)\, \dot\eta
  \, d\tau- \int_{a_k}^{b_k} H_\gamma(\eta,
 U'(\eta)) \, d\tau = U(\eta(b_k))-U(\eta(a_k)) + c \, (b_k-a_k)
\end{eqnarray*}
for any $k$. The indices $k$ for which the  curve $\eta$ restricted to
 $[a_k,b_k]$ is noncyclic are finitely many because the length of $\eta$ is finite,
  when we sum
 the contributions $U(\eta(b_k))-U(\eta(a_k))$ over this family of indices we get
  $U(\eta(t_1) -U(\eta(t_0))$,  where
\begin{eqnarray*}
  t_0 &=& \inf\{t \in [a,b] \mid \eta(t)\in E\}, \\
  t_1 &=& \sup\{t \in [a,b] \mid \eta(t)\in E\}
\end{eqnarray*}
If instead $\eta$ restricted to $[a_k,b_k]$ is cyclic then clearly $U(\eta(b_k))-U(\eta(a_k))=0$.
Altogether  we have
\begin{eqnarray}
  \int_E \bar L_\gamma(\eta,\dot\eta) \, d\tau &=&
  \sum_k \int_{a_k}^{b_k} \bar L_\gamma(\eta,\dot\eta) \, d\tau\geq \sum_k (U(\eta(b_k))-U(\eta(a_k))) + c \, (b_k-a_k) \label{cruxtris2}\\ &=&
    U(\eta(t_1))-U(\eta(t_0)) + c \, |E|. \nonumber
\end{eqnarray}
Taking into account that
\[\dot \eta(t)= 0 \qquad\hbox{for a.e. $t \in F$}\]
we further get
\begin{equation}\label{cruxtris3}
  \int_F \bar L_\gamma(\eta,\dot\eta) \, d\tau = c \, |F|.
\end{equation}
If $a,\,b\in E$, then  $t_0=a$, $t_1=b$ and by the assumption
\begin{equation}\label{cruxtris4}
 U(\eta(t_0))=U(\eta(a))=U(\eta(b))=U(\eta(t_1)).
\end{equation}
If instead $a,\,b\in F$ then $U(\eta(t_0))=U(\eta(a))$ and $U(\eta(t_1))=U(\eta(b))$ and
we again obtain \eqref{cruxtris4}.
By \eqref{cruxtris2}, \eqref{cruxtris3}, \eqref{cruxtris4} we conclude that
\[\int_a^b \bar L_\gamma(\eta,\dot\eta) \, d\tau \geq c \, (b-a),\]
as was asserted.
\endproof

\smallskip

\begin{rem}\label{postpost} If we have in the statement of Proposition \ref{postprecruxtris}
the additional information  that  $\eta([a,b]) \cap \{0\} = \emptyset$
 and set
 \begin{equation}\label{postpost1}
  F(t) = \left\{ \begin{array}{cc}
                    L_\gamma(\eta(t),\dot\eta(t)) & \quad\hbox{for $\eta(t) \neq 1$}  \\
                    c_{\gamma(1)} & \quad\hbox{for $\eta(t)=1$} \\
                  \end{array} \right .
 \end{equation}
the same argument in the proof of the proposition
allows showing
\[\int_a^b  F(\tau) \, d\tau \geq
  c_{\gamma(1)} \, (b-a) .\]
  \end{rem}

\bigskip

\section{The action functional} \label{act}

\subsection{The Lagrangian on \texorpdfstring{$\Gamma$}{Γ}} \label{act1}
 We assume that it is given  a flux limiter
$c_x$ for any $x \in \Vbf$. We  set for $(x,q) \in \Vbf \times \Rds$
 \[ \Ecal(x,q) =  \{ \gamma \in \Gamma_x \mid \;
 \hbox{$q$ parallel to $\dot\gamma(1)$}\},\]
 note that $\Ecal(x,q) \neq \emptyset$ if and only if $(x,q) \in T\Gamma$.
In the definition of  the Lagrangian $L: T\Gamma \to \Rds$ we distinguish three cases:
\begin{itemize}
  \item [--] $x \in \Gamma \setminus \Vbf$, $x \in \gamma((0,1))$, then
  \begin{equation}\label{lagrange}
 L(x,q) = L_\gamma \left (\gamma^{-1}(x),\frac{q \cdot \dot\gamma(\gamma^{-1}(x))}
 {|\dot \gamma(\gamma^{-1}(x))|^2} \right ) ;
  \end{equation}
  \item[--] $x \in \Vbf$, $q \neq 0$, then
  \begin{equation}\label{lagrangebis}
   L(x,q) = \min_{\gamma \in \Ecal(x,q)} L_{\gamma} \left (1, \frac{q \cdot \dot\gamma
 (1)}{|\dot \gamma(1)|^2} \right );
  \end{equation}
  \item[--] $x \in \Vbf$, $q = 0$, then
  \begin{equation}\label{lagrangetris}
    L(x,0)= c_x .
  \end{equation}
  \end{itemize}

\smallskip

Note that  formula \eqref{lagrangebis} is more involved
because there is a problem to take into account,
  namely different arcs ending at $x$ could have parallel tangent vectors, in this case
we should have
\[q= \lambda_1 \, \dot \gamma_1(1) = \lambda_2 \, \dot\gamma_2(1) \qquad
\hbox{for arcs $\gamma_1 \neq \gamma_2$, scalars $\lambda_1,\, \lambda_2$.}\]
As already pointed out in the Introduction, \eqref{lagrangetris} provides
the link between flux limiter and Lagrangian, and consequently between flux limiter
 and representation  formula.

\smallskip

\begin{prop} The Lagrangian $L: T\Gamma \to \Rds$ is lower semicontinuous.
\end{prop}
\proof It is clearly continuous in $T\Gamma \cap \big ( (\Gamma \setminus \Vbf) \times \Rds^N \big )$. Now consider  $x \in \Vbf$, $q \neq 0$, and $(x_n,q_n) \in T\Gamma$, $(x_n,q_n) \to (x,q)$. Taking into account that the arcs are finitely many, we can assume, up  to extracting  a subsequence, that $x_n = \gamma(s_n)$ for an arc $\gamma \in \Gamma_x$ with $s_n \to 1$, $q_n=\mu_n \, \dot \gamma(s_n)$, for some scalar $\mu_n$, and $L(x_n,q_n)= L_\gamma(s_n,\mu_n)$. Hence $\mu_n \to \mu$, $q= \mu \,  \dot\gamma(1)$ and
\[L(x_n,q_n)= L_\gamma(s_n,\mu_n) \to L_\gamma(1,\mu) \geq L(x,q).\]
Now assume $q=0$ and consider $(x_n,q_n) \to (x,0)$, if   $q_n \neq 0$ or $x_n \neq x$,
we can argue as above, bearing in mind Lemma \ref{flusso}.
It is left the case $(x_n,q_n) \equiv (x_0,0)$, which is trivial.
\endproof

\smallskip

It is convenient to extend $L$ on the whole of $\Rds^N \times \Rds^N$, keeping it lower semicontinuous, setting
\[L(x,q)= + \infty  \qquad\hbox{if $(x,q) \not\in T\Gamma$.}\]

It is apparent that $L$ is bounded from below. We set

\begin{equation}\label{ml}
  m_L = \min_{(x,q) \in T\Gamma} L(x,q).
\end{equation}

From the superlinearity of the $L_\gamma$'s we finally derive that  there exists
a function $\theta: \Rds^+ \to \Rds$ with $\lim_{t \to + \infty} \frac {\theta(t)}t= + \infty$
such that
\begin{equation}\label{theta}
  L(x,q) \geq \theta(|q|) \qquad\hbox {for any $(x,q) \in T\Gamma$.}
\end{equation}

\medskip
\subsection{Admissible curves}
Given a curve $\xi: [0,T] \to \Rds^N$, we define the corresponding action functional as
\[\int_0^T L(\xi,\dot \xi) \, d\tau.\]

\smallskip

Note that for any pair of points $x$ and $y$ of $\Gamma$ and $T>0$ there are curves linking $x$ to $y$ in the time $T$ with finite action functional. It is enough for that to take a geodesic linking $x$ to $y$, which does exist since the network is connected, and to change the parametrization in order to define it in $[0,T]$.

\smallskip

The following result is a consequence of Proposition \ref{postprecruxtris} and
Remark \ref{postpost}.
\begin{prop}\label{precrux} Let $\xi$
be a curve defined in $[a,b]$, $x \in \Vbf$.
Assume  that   $\xi(a)= \xi(b)$,  and $\xi([a,b]) \subset \gamma([0,1])$ for some
 $\gamma \in \Gamma_x$. Assume further that
  $\xi([a,b]) \cap (\Vbf \setminus \{x\}) = \emptyset$, then
\[\int_a^b L(\xi,\dot\xi) \, d\tau \geq c_x \, (b-a).\]
\end{prop}
\proof
We set
\begin{align*}
 E=&\,\{t \in [a,b] \mid \xi(t) \neq x\},\\
 F=&\,\{t \in [a,b] \mid \xi(t) = x\}
\end{align*}
and define
\[\eta(t)= \gamma^{-1}\circ \xi(t)  \qquad\hbox{for $t \in [a,b]$}.\]
Note that $\eta$ satisfies the assumptions of Proposition \ref{postprecruxtris} plus
$\eta([a,b]) \cap \{0\} = \emptyset$.
Taking into account the definition of the Lagrangian $L$ and Lemma \ref{deri}, we have that
\begin{equation}\label{precrux1}
 L(\xi(t),\dot\xi(t))= L_\gamma(\eta(t), \dot\eta(t)) \qquad\hbox{for a.e. $t \in E$,}
\end{equation}
in addition
\begin{eqnarray}
  \dot\xi(t) &=&0 \qquad \hbox{for  a.e. $t \in F$} \label{precrux2} \\
  L(x,0) &=& c_x.  \label{precrux3}
\end{eqnarray}
If we define $F(t)$, for $t \in [a,b]$, as in \eqref{postpost1}, we get by \eqref{precrux1},
\eqref{precrux2}, \eqref{precrux3}
\[F(t)= L(\xi(t),\dot\xi(t)) \qquad\hbox{ for a.e. $t$}\]
and the assertion is then a consequence of Proposition \ref{postprecruxtris},
Remark \ref{postpost}.
\endproof

\smallskip

\begin{defin}\label{admissible}
We say that a curve $\xi:[0,T] \to \Gamma$ is {\em admissible} if
 there exists a finite partition $\{t_1,\dotsc, t_m\}$ of the interval $[0,T]$ with
\begin{itemize}
  \item[--] $t_1=0$, $t_m=T$, $\xi(t_i) \in \Vbf$ for $i=2, \dotsc, m-1$ if $m >2$;
  \item[--] for any $i$, either $\xi((t_i,t_{i+1})) \cap \Vbf = \emptyset$  and $\xi(t_i) \neq \xi(t_{i+1})$ or $\xi(t) \equiv x \in \Vbf$ for $t \in (t_i,t_{i+1})$.
\end{itemize}
\end{defin}
\smallskip
\begin{lem}\label{admi} Let $\xi: [0,T] \to \Gamma$ be a
 curve. Assume that for any nondegenerate interval $[a,b]$ with $\xi(a) ,\, \xi(b) \in \Vbf$,
$\xi((a,b)) \cap \Vbf = \emptyset$, one
has $\xi(a) \neq \xi(b)$ then $\xi$ is admissible.
\end{lem}
\proof
We will show that, under the assumption in the statement, there is a finite partition of $[0,T]$ satisfying Definition \ref{admissible}. If $\xi((0,T))\cap\Vbf=\emptyset$, then $\{0,T\}$ is such a partition. Otherwise let us denote by $r>0$ the minimum of the geodesic distance $d_\Gamma(x,y)$ between distinct vertices of the network. Being $\xi$ absolutely continuous, there exists a $\delta>0$ such that
\[
\int_E\left|\dot\xi(\tau)\right|d\tau<r,\qquad\text{if }|E|\le\delta.
\]
Let us partition $[0,T]$ into a finite number of closed intervals $I_1,\dotsc,I_n$ of length less than $\delta$. Then $\xi(I_k)$ can contain at most one vertex $x_k$, for each $k\in\{1,\dotsc,n\}$. Let us set
\[
E_k:=\{t\in I_k:\xi(t)\ne x_k\}
\]
and pick $k$ such that $E_k\ne I_k$. Then $E_k$ is open in $I_k=[a_k,b_k]$. Let us assume there exists a connected component $(a,b)$ of $E_k$ with $a_k<a$ and $b<b_k$. Then $\xi(a)=\xi(b)=x_k$, hence by assumption $\xi((a,b))=\{x_k\}$, in contradiction with the definition of $E_k$. This shows that, if not empty,
\begin{equation}\label{eq:admi1}
E_k=[a_k,\alpha_k)\cup(\beta_k,b_k],\qquad\text{for some }a_k\le\alpha_k\le\beta_k\le b_k.
\end{equation}
Finally \eqref{eq:admi1} implies that
\[
\{t\in[0,T]:\xi(t)\notin\Vbf\}=\bigcup_{k=1}^nE_k
\]
is the disjoint union of a finite number of intervals whose extreme points, possibly adding $0$ and $T$, form a partition of $[0,T]$ satisfying Definition \ref{admissible}.
\endproof

\smallskip

\begin{prop}\label{crux} Given $T >0$, and a curve $\xi$ defined in $[0,T]$,
we can find an admissible curve $\zeta$ defined in the same interval with $\xi(0)=\zeta(0)$, $\xi(T)=\zeta(T)$ and
\begin{equation}\label{eq:crux.1}
\int_0^T L(\xi,\dot \xi) \, d\tau \geq \int_0^T L(\zeta,\dot \zeta) \, d\tau
\end{equation}
\end{prop}
\proof
We can assume that $\xi$ is supported in $\Gamma$, otherwise the action
is infinite and there is nothing to prove. If $\xi$
is admissible then, again,  we are done. If not, there is,
according to Lemma \ref{admi}, a collection of  nondegenerate  intervals
$[a_i,b_i] \subset [0,T]$ with disjoint interiors such that
$\xi(a_i)=\xi(b_i) \in \Vbf$ and $\xi((a_i,b_i)) \cap \Vbf = \emptyset$.
Note that they are not more than countably many because $\Rds$ is separable.
We fix $i=1$ and
 find an arc $\gamma \in \Gamma_x$ with $\xi(a_1,b_1) \subset \gamma((0,1))$.
 We define
 \[\zeta(t)= \left \{ \begin{array}{cc}
               x & \quad t \in [a_1,b_1] \\
               \xi(t) & \quad t \in [0,T] \setminus [a_1,b_1]
             \end{array} \right . \]
and  get  from Proposition \ref{precrux} and the very definition of $L$
\begin{eqnarray*}
 \int_0^T L(\zeta, \dot \zeta) \, d\tau  &=& \int_0^{a_1} L(\xi, \dot \xi) \, d\tau
 + c_x \, (b_1-a_1) + \int_{b_1}^T L(\xi, \dot \xi) \, d\tau  \leq \int_0^T L(\xi, \dot \xi) \, d\tau.
\end{eqnarray*}
We can repeat the above procedure for any other $[a_i,b_i]$ for $i >1$. Note that $[0,T] \setminus \cup_i [a_i,b_i]$  is made up by   intervals which  have disjoint interiors  and endpoints at two different vertices, except possibly the interval starting at $\xi(0)$ and the one ending at $\xi(T)$.   These intervals must be finitely many  because $\xi$ is absolutely continuous and so possess finite length.   Consequently the $[a_i,b_i]$'s can be glued together in a finite number of disjoint intervals corresponding to steady states at a given vertex.
We have  therefore constructed in this way an admissible curve satisfying
the statement.
\endproof
\endproof
\smallskip

Given $\xi:[0,T] \to \Gamma$, we will denote in what follows by $\xi_{\mathrm {ad}}$
 the admissible curve
obtained from $\xi$ through the procedure detailed  in the above proof.

\smallskip

\begin{rem}\label{remcrux} Looking back at
  the construction performed in the proof of the Proposition \ref{crux},
  it is apparent that
\begin{itemize}
  \item[{\bf (i)}]  $\xi_{\mathrm {ad}}(t) \in \Gamma \setminus \Vbf \, \Rightarrow \,
  \xi_{\mathrm {ad}}(t) =\xi(t)$;
  \item[{\bf (ii)}]  $\xi(t) \in  \Vbf \, \Rightarrow \,
  \xi_{\mathrm {ad}}(t) =\xi(t)$;
  \item[{\bf (iii)}] if $\{t_i\}$ is  an admissible partition of
  $[0,T]$ associated to $\xi_{\mathrm {ad}}$,
      then $\xi_{\mathrm {ad}}(t_i)=\xi(t_i)$.
\end{itemize}

\end{rem}

\smallskip

\begin{prop}\label{postcrux}
Let $\xi:[a,b]\to\Gamma$ be an admissible curve and assume that there is an $x\in\Vbf$ such that
\begin{equation}\label{eq:postcrux.1}
\xi([a,b])\cap\Vbf=\{x\}
\end{equation}
or
\begin{equation}\label{eq:postcrux.2}
\xi([a,b])\subseteq\gamma((0,1]),\qquad\text{for some }\gamma\in\Gamma_x.
\end{equation}
Then there exists a positive constant $\ell$ such that
\[
\int_a^bL(\xi,\dot\xi)\,d\tau\ge c_x\,(b-a)-\ell\,(d_\Gamma(x,\xi(a))\vee d_\Gamma(x,\xi(b))).
\]
\end{prop}
\proof
We recall from Lemma \ref{lippi} that $\frac1m$, with $m$ defined as in \eqref{lippi1}, is a Lipschitz constants of $\gamma^{-1}:\gamma([0,1])\to[0,1]$ for all arcs $\gamma$. We further set
\[
M=\max_{\gamma\in\Ecal}\max_{s\in[0,1]}L(\gamma(s),\dot\gamma(s))\vee1.
\]
We will show that the inequality in the statement holds true with
\[
\ell:=2\left(\frac{M-\min\limits_{x\in\Vbf}(c_x\wedge0)}m\right).
\]
We start looking into the case \eqref{eq:postcrux.2}. We assume, to fix ideas, that
\[
\alpha:=\gamma^{-1}\circ\xi(a)\le\gamma^{-1}\circ\xi(b)=:\beta.
\]
We then consider the curve $\zeta$ obtained via concatenation of $\xi$ and $\widetilde\gamma$ restricted to $[1-\beta,1-\alpha]$, we have
\[
1-\alpha\le\frac dm\qquad\text{and }\qquad1-\beta\le\frac dm,
\]
where $d=d_\Gamma(x,\xi(a))\vee d_\Gamma(x,\xi(b))$. We can apply to $\zeta$ Proposition \ref{precrux} to get
\[
\int_a^{b+\beta-\alpha}L(\zeta,\dot\zeta)\,d\tau\ge c_x\,(b-a+\beta-\alpha).
\]
We further have
\[
\int_b^{b+\beta-\alpha}L(\zeta,\dot\zeta)\,d\tau=\int_{1-\beta}^{1-\alpha}L(\widetilde\gamma,\dot{\widetilde \gamma})\,d\tau\le M\,(\beta-\alpha)\le M\,(1-\alpha)\le\frac Mm\,d
\]
and accordingly, if $c_x\ge0$,
\[
\int_a^bL(\xi,\dot\xi)\,d\tau\ge c_x\,(b-a)+c_x\,(\beta-\alpha)-\frac Mm\,d\ge c_x\,(b-a)-\frac Mm\,d=c_x\,(b-a)-\frac\ell2\,d
\]
or, if $c_x<0$,
\[
\int_a^bL(\xi,\dot\xi)\,d\tau\ge c_x\,(b-a)+c_x\,(1-\alpha)-\frac Mm\,d\ge c_x\,(b-a)-\left(\frac{M-c_x}m\right)d\ge c_x\,(b-a)-\frac\ell2\,d.
\]
In both cases we then get
\begin{equation}\label{eq:postcrux1}
\int_a^bL(\xi,\dot\xi)\,d\tau\ge c_x\,(b-a)-\frac\ell2\,d.
\end{equation}
The case where $\gamma^{-1}\circ\xi(a)\ge\gamma^{-1}\circ\xi(b)$ can be handled similarly.\\
Now let consider the case \eqref{eq:postcrux.1}. By the assumption, $\xi$ visits exactly two different arcs in the time interval $[a,b]$. It starts in one arc at time $a$, then intersects the vertex $x$, and possibly stops on it for a positive time, and finally enters the second arc, where it dwells until the time $b$. It cannot go back to $x$, otherwise the admissibility condition should be violated. To fix ideas, assume that
\[
\xi([a,b])\cap\gamma_1((0,1))\ne\emptyset,\qquad\xi([a,b])\cap\gamma_2((0,1))\ne\emptyset,
\]
where $\gamma_1$ and $\gamma_2$ are different arcs, both belonging to $\Gamma_x$, and that $\xi(a)\in\gamma_1((0,1))$, $\xi(b)\in\gamma_2((0,1))$. We thus have three disjoint intervals $I_1$, $I_2$ and $I_3$ covering $[a,b]$ and such that
\[
\xi(I_1)\subseteq\gamma_1((0,1)),\qquad\xi(I_2)\subseteq\gamma_2((0,1)),\qquad\xi(I_3)\equiv x.
\]
By definition
\[
\int_{I_3}L(\xi,\dot\xi)\,d\tau=c_x|I_3|
\]
and it follows from \eqref{eq:postcrux1} that
\[
\int_{I_1}L(\xi,\dot\xi)\,d\tau+\int_{I_2}L(\xi,\dot\xi)\,d\tau\ge c_x(|I_1|+|I_2|)-\ell(d_\Gamma(x,\xi(a))\vee d_\Gamma(x,\xi(b))).
\]
This concludes the proof.
\endproof

\bigskip

\section{Minimal action functional} \label{min}

We proceed proving  a semicontinuity property for the action functional with respect to  the
weak topology of $H^{1,1}([0,T])$, the Sobolev space of the absolutely continuous curves with the norm
\[
\|\xi\|:=\int_0^T(|\xi|+|\dot\xi|)d\tau.
\]
The difficulty here is that the
 Lagrangian $L$ is not convex in the velocity variable at the vertices.
A deeper analysis using Proposition \ref{postcrux} is therefore needed.

\smallskip
\begin{prop}\label{debole}  Let $\xi_n$ be a sequence of admissible curves, defined in $[0,T]$,
weakly converging to $\xi$ in $H^{1,1}$, then
 \[\liminf_n \int_0^T  L(\xi_n,\dot \xi_n) \, d\tau \geq
 \int_0^T  L(\xi_{\mathrm {ad}},\dot \xi_{\mathrm {ad}}) \, d\tau.\]
\end{prop}
\proof We fix $\eps >0$. Taking into account that the $\dot\xi_n$'s  are uniformly integrable
by the convergence assumption on $\xi_n$,
we can choose  $\delta \in (0,\eps)$ such that
\begin{equation}\label{debole0}
 |E| < 2 \, \delta
\, \Rightarrow \, \left \{\begin{array}{cc}
 \int_E  |\dot\xi_n| \, d\tau < \eps &  \quad \forall \, n \\
 \int_E |L(\xi_{\mathrm {ad}}, \dot \xi_{\mathrm {ad}})| \, d\tau < \eps&
 \end{array} \right .
\end{equation}
Since $\xi_{\mathrm {ad}}$ is admissible,
$\xi_{\mathrm { ad}}^{-1}(\Vbf)$  is made up  by finitely many
compact intervals, some of them reduced to a singleton,
and some other nondegenerate,  where $\xi_{\mathrm {ad}}$ is constant, say
\[ \xi_{\mathrm {ad}}^{-1}(\Vbf)= \cup_{k=1}^m J_k.\]
We  denote by $x_k$,
$k=1, \dotsc, m$, the element of $\Vbf$  univocally determined by the condition
\[x_k= \xi_{\mathrm {ad}}(t) \qquad\hbox{for some $t \in J_k$.}\]
Accordingly
\[A_\delta= \{t \in [0,T] \mid d(t, \xi_{\mathrm {ad}}^{-1}(\Vbf)) \leq \delta\}\]
is a finite union of disjoint compact intervals denoted by $J'_1, \dotsc, J'_m$,
 with $[a_k,b_k] =:J'_k \supset J_k$,
provided that $\eps$, and so $\delta$, is sufficiently small.
Since $\xi_n$ uniformly converges to $\xi$, we deduce, taking into account
Remark \ref{remcrux} {\bf (ii)}, that  the $\xi_n$'s
restricted to any $J'_k$ satisfy the assumption of Proposition
 \ref{postcrux} for $n$ large enough, with $x=x_k$. In addition
 \begin{eqnarray*}
   d_\Gamma(\xi_n(a_k),x_k) &\leq & \int_{a_k}^{a_k+\delta} |\dot\xi_n| \, d\tau \leq \eps \\
   d_\Gamma(\xi_n(b_k),x_k) &\leq & \int_{b_k- \delta}^{b_k} |\dot\xi_n| \, d\tau \leq \eps
 \end{eqnarray*}
 for any such $n$, and any $k$, by \eqref{debole0}.
 We consequently  have by Proposition \ref{postcrux}
 \begin{equation}\label{debole1}
  \int_{J'_k} L(\xi_n, \dot\xi_n) \, d\tau \geq c_{x_k} \, |J'_k|- \ell \, \eps.
 \end{equation}
 On the other side, by the very definition of $J_k$ and \eqref{debole0}, we have
 \begin{equation}\label{debole2}
  \int_{J'_k} L(\xi_{\mathrm {ad}}, \dot\xi_{\mathrm {ad}}) \, d\tau =
  \int_{J_k} L(\xi_{\mathrm {ad}}, \dot\xi_{\mathrm {ad}}) \, d\tau + \int_{J'_k \setminus J_k}
  L(\xi_{\mathrm {ad}}, \dot\xi_{\mathrm {ad}}) \, d\tau \nonumber \leq  c_{x_k} \,|J_k|  + \eps  \leq
    c_{x_k} \,|J'_k| + \eps.
 \end{equation}
 By combining \eqref{debole1}, \eqref{debole2} we get
 \begin{equation}\label{debole3}
  \int_{J'_k} L(\xi_n, \dot\xi_n) \, d\tau
 \geq \int_{J'_k} L(\xi_{\mathrm {ad}}, \dot\xi_{\mathrm {ad}}) \, d\tau - (\ell +1) \, \eps,
 \end{equation}
for any $k=1, \dotsc, m$, $n$ sufficiently large. The complement in $[0,T]$ of
$\overline{A_\delta}$  is a finite union of open (in $[0,T]$)  intervals denoted by
$I_1, \dotsc, I_l$, where $\xi_{\mathrm {ad}}$ and $\xi$ coincide by Remark \ref{remcrux}.
Note that on the intervals $I_r$, the Lagrangian satisfies the usual conditions of
lower semicontinuity and convexity which imply  that the action functional is
sequentially  lower semicontinuous in the weak topology of $H^{1,1}$, see Theorem 3.6 in \cite{Buttazzo},
 so that we get
\begin{equation}\label{debole4}
  \liminf_{n \to +\infty} \int_{I_r} L(\xi_n,\dot\xi_n) \, d\tau
  \geq \int_{I_r} L(\xi_{\mathrm {ad}},\dot\xi_{\mathrm {ad}}) \, d\tau
  \quad\hbox{for any $r= 1, \dotsc, l$.}
\end{equation}
By combining  \eqref{debole3}, \eqref{debole4}, we further get
\[\liminf_{n \to +\infty} \int_0^T L(\xi_n,\dot\xi_n) \, d\tau \geq
\int_0^T L(\xi_{\mathrm {ad}},\dot\xi_{\mathrm {ad}}) \, d\tau - m \, (\ell+1)\,\eps,\]
which  gives the assertion since $\eps$ has been arbitrarily chosen.
\endproof

\smallskip
\begin{theo}\label{funda} Given $x$ and $y$ in $\Gamma$ and $T>0$, there is an
admissible curve defined in $[0,T]$ and linking $x$ to $y$ which minimizes  the action.
\end{theo}
\proof Let $\xi_n$ be a minimizing sequence that we can assume, according
to Proposition \ref{crux},
made up by admissible curves. Since the minimal action in $[0,T]$ between $x$, $y$
is finite, say equal to $M$, we find
\[M +1 \geq \int_0^T L(\xi_n,\dot \xi_n) \, dt \geq \int_0^T \theta(|\dot \xi_n|) \, dt,\]
for $n$ large, where $\theta$ has been introduced in \eqref{theta}.
Since  the $\dot\xi_n$'s  are  bounded in $L^1$,
we can use the Dunford--Pettis compactness criterion, see Theorem 4.7.18 in \cite{Bogachev07}, to get that the $\xi_n$ converges
weakly in $H^{1,1}$, up to subsequences, to a limit curve $\xi$.
By Proposition \ref{debole}
\[\liminf_n \int_0^T L(\xi_n,\dot\xi_n) \, d\tau \geq \int_0^T L(\xi_{\mathrm {ad}},\dot\xi_{\mathrm {ad}}) \, d\tau.\]
This shows  that $\xi_{\mathrm {ad}}$ is the sought minimizer.
\endproof

\smallskip
Next result is relevant by itself.

\smallskip
\begin{theo}\label{lippo} Any admissible curve minimizing the action between two elements $z_1$,
$z_2$ of $\Gamma$ in a given time $T$ is Lipschitz continuous.
\end{theo}
\begin{proof} Let $\xi:[0,T] \to \Gamma$ be an admissible minimizing curve and $\{t_1, \cdots, t_m\}$ the corresponding finite partition of $[0,T]$. We will show that the restriction of $\xi$ to any interval $[t_i,t_{i+1}]$, $i= 1,\cdots, m-1$ is Lipschitz continuous. If such a restriction is constant there is clearly nothing to prove. We  therefore focus on an interval $[t_i,t_{i+1}]$ with
\[\xi((t_i,t_{i+1})) \subset \gamma((0,1)) \qquad\hbox{for a suitable arc $\gamma$.}\]
We consequently have
\[\int_{t_i}^{t_{i+1}} L (\xi, \dot \xi) \, dt = \int_{t_i}^{t_{i+1}} L_\gamma (\eta, \dot \eta) \, dt,\]
where $\eta=\gamma^{-1}\circ \xi$, see Lemma \ref{deri}.  To fit the setting of Section \ref{covidcovid} we can extend $L_\gamma$ in the whole of $\Rds \times \Rds$ defining for instance
\[L_\gamma(s,\lambda)= \left \{\begin{array}{c}
                      L_\gamma(0,\lambda) \quad\hbox{if $s <0$} \\
                      L_\gamma(1,\lambda) \quad\hbox{if $s >1$}
                    \end{array} \right .\]
By the optimality properties of $\xi$, it is clear that $\eta$ minimizes the action induced by $L_\gamma$ among the curves $\zeta$ with $\zeta(t_i)=\eta(t_i)$, $\zeta(t_{i+1})=\eta(t_{i+1})$, $\spt \zeta \subset \spt \eta$. We can therefore invoke Theorem \ref{covid} to deduce that $\eta$, and consequently $\xi$, is Lipschitz continuous in $[t_1,t_2]$.
\end{proof}

\smallskip

We proceed proving:

\begin{theo} The minimal action functional
\[(x,t,y,r) \mapsto \min \left \{\int_0^{r-t} L(\xi,\dot\xi) \, d\tau \mid \xi\;\hbox{AC with}\;
\xi(0)=x,\, \xi(r-t)=y \right\}\]
is continuous for $x$, $y$ varying in $\Gamma$, $t$, $r$ in $\Rds^+$ with $r >t$.
\end{theo}
\proof We denote by $S(\cdot,\cdot,\cdot,\cdot)$ the minimal action
functional. We first show that it is lower semicontinuous.
We fix $(x_0,t_0,y_0,r_0)$ with $r_0 > t_0$, and consider
$(x_n,t_n,y_n,r_n)  \to (x_0,t_0,y_0,r_0)$ with
\[S(x_n,t_n,y_n,r_n) \to a \qquad\hbox{for some $a \in \Rds \cup \{+\infty\}$.}\]
If $a= +\infty$ there is nothing to prove, we can therefore assume that
\begin{equation}\label{conticonti01}
 \int_0^{r_n-t_n} L(\xi_n,\dot\xi_n) \, d\tau  \quad\hbox{is bounded from above,}
\end{equation}
where   $\xi_n$  denotes a sequence of admissible
curves defined
in $[0,r_n-t_n]$ and realizing the minimal action between $x_n$ and $y_n$.
Given $\eps >0$, we set $T_\eps=r_0-t_0+\eps$, so that $r_n-t_n < T_\eps$ for $n$ large,
 and we extend all such $\xi_n$'s in $[0,T_\eps]$ defining
\[\xi_n(t) = y_n  \qquad\hbox{for $t \in (r_n-t_n,T_\eps]$}.\]
Arguing as in the proof of Theorem \ref{funda} and taking into account \eqref{conticonti01},
 we see that the $\xi_n$'s
 weakly converge in $H^{1,1}((0,T_\eps))$ to a curve $\xi$ with $\xi(0)=x$, $\xi(r-s)=y$.
 By Proposition \ref{debole}, we have
 \[\liminf_n \int_0^{T_\eps} L(\xi_n,\dot\xi_n) \, d\tau
 \geq \int_0^{T_\eps} L(\xi_{\mathrm {ad}},\dot\xi_{\mathrm {ad}}) \, d\tau.\]
 We further have
 \begin{eqnarray*}
  \int_0^{T_\eps} L(\xi_n,\dot\xi_n) \, d\tau &=& S(x_n,t_n,y_n,r_n) +
  \int_{r_n-t_n}^{T_\eps} L(y_n,0) \, d\tau \\
    &\leq& S(x_n,t_n,y_n,r_n) +
  \widetilde M \, (r_0-t_0+\eps -r_n+t_n),
 \end{eqnarray*}
 where $\widetilde M= \max_{x \in \Gamma} |L(x,0)|$,  and consequently
 \[\int_0^{T_\eps} L(\xi_n,\dot\xi_n) \, d\tau \leq S(x_n,t_n,y_n,r_n)
 + 2 \, \widetilde M \, \eps\]
 for $n$ suitably large. We deduce that
 \begin{eqnarray*}
 a &\geq& \liminf_n \int_0^{T_\eps} L(\xi_n,\dot\xi_n) \,
 d\tau - 2 \, \widetilde M \, \eps \geq  \int_0^{T_\eps} L(\xi_{\mathrm {ad}},\dot\xi_{\mathrm {ad}}) \, d\tau
    - 2 \, \widetilde M \, \eps\\
     &\geq& S(x_0,t_0,y_0,r_0) + (m_L - 2 \, \widetilde M) \, \eps,
 \end{eqnarray*}
 where $m_L$ is defined as in \eqref{ml}.
 This shows the claimed lower semicontinuity.

 We proceed considering again $(x_n,t_n,y_n,r_n)  \to (x_0,t_0,y_0,r_0)$ with
\[S(x_n,t_n,y_n,r_n) \to a    \qquad\hbox{for some $a $,} \]
and denote by $\xi:[0,r_0-t_0] \to \Gamma$
 an admissible curve realizing the minimal
action between  $x_0$ and $y_0$.
We further set $\tau_k = \frac 1k$, $\tau'_k = (r_0-t_0) - \frac 1k$,
the velocity $|\dot\xi|$ being bounded by Theorem \ref{lippo}, we  have
\begin{eqnarray}
  \int_0^{\tau_k} L(\xi,\dot\xi) \, d\tau &=& \mathrm O(1/k) \label{conticonti00}\\
  \int_{\tau'_k}^{r_0-t_0} L(\xi,\dot\xi) \, d\tau &=& \mathrm O(1/k), \label{conticonti0}
\end{eqnarray}
where $\mathrm O(\cdot)$ stands for the Landau symbol, and
\[d_\Gamma(\xi(\tau_k),x_0) = \mathrm O(1/k), \quad d_\Gamma(\xi(\tau'_k),y_0) = \mathrm O(1/k).\]
We  select $n_k$
such that
\[|r_{n_k} - t_{n_k} -r +s| =\mathrm O(1/k^2), \quad d_\Gamma(x_{n_k},x_0)= \mathrm O(1/k^2)
, \quad d_\Gamma(y_{n_k},y_0)=\mathrm O(1/k^2).\]
We therefore have
\[ d_\Gamma(x_{n_k},\xi(\tau_k))= \mathrm O (1/k) \qquad\hbox{and}
\qquad d_\Gamma(y_{n_k},\xi(\tau'_k))= \mathrm O (1/k).\]
Let  $\zeta_k$ be a geodetic in $\Gamma$ linking $\xi(\tau_k)$ and $x_{n_k}$
reparametrized in order to be defined in the interval $[0,\tau_k]$. We have
\begin{equation}\label{conticonti1}
 \int_0^{\tau_k} L(\zeta_k,\dot\zeta_k) \, d\tau = \mathrm O(1/k).
\end{equation}
Let  $\zeta'_k$ be a geodetic linking $\xi(\tau'_k)$ and $y_{n_k}$
reparametrized in order to be defined in the interval $[\tau'_k,r_{n_k} -t_{n_k}]$. We have
\begin{equation}\label{conticonti2}
 \int_{\tau'_k}^{r_{n_k}-t_{n_k}} L(\zeta'_k,\dot\zeta'_k) \, d\tau = \mathrm O(1/k).
\end{equation}
The curve
\[\xi_k(\tau)= \left \{ \begin{array}{cc}
           \zeta_k(\tau) &  \quad \tau \in [0,\tau_k] \\
           \xi(\tau) & \quad \tau \in (\tau_k,\tau'_k) \\
           \zeta'_k(\tau) & \quad \tau \in [\tau'_k,r_{n_k} -t_{n_k}]
         \end{array} \right . \]
links $x_{n_k}$ to  $y_{n_k}$ in the time interval $[0,r_{n_k} -t_{n_k}]$  for $k$ large.
 We have by \eqref{conticonti1}, \eqref{conticonti2}, \eqref{conticonti00},
 \eqref{conticonti0}
\begin{eqnarray*}
 \int_0^{r_{n_k} -t_{n_k}} L(\xi_k,\dot\xi_k) \, d\tau &=&
 \int_{\tau_k}^{\tau'_k} L(\xi,\dot\xi) \, d\tau+ \mathrm O(1/k) \\
   &=& \int_0^{r-s} L(\xi,\dot\xi) \, d\tau+ \mathrm O(1/k).
\end{eqnarray*}
This implies
\begin{eqnarray*}
  a &=& \lim_k S(x_{n_k},t_{n_k},y_{n_k},r_{n_k})  \\
  &\leq&  \lim_k \int_0^{r_{n_k} -t_{n_k}} L(\xi_k,\dot\xi_k) \, d\tau = S(x_0,t_0,y_0,r_0),
\end{eqnarray*}
which shows the upper semicontinuity of the minimal action.
\endproof

\bigskip

\section{The time--dependent equation on \texorpdfstring{$\Gamma$}{Γ}}  \label{mai}

For any given arc $\gamma$, we consider  the time--dependent
equation
\begin{equation}\label{HJg}  \tag{HJ$_\gamma$}
    U_t + H_\gamma(s,U')=0 \qquad\hbox{in $ \Qcal$.}
\end{equation}
We are interested in finding  a continuous function $v: \Gamma \times \Rds^+ \to \Rds$
 such that $v \circ \gamma$ solves \eqref{HJg} in $\Qcal$, for any $\gamma$, taking
 into account, in the sense we are going to specify, a flux limiter $c_x$
  at any vertex. We denote by (HJ$\Gamma$) the problem as a whole.
\smallskip

The definition of (sub / super) solution to (HJ$\Gamma$) is as follows:

\smallskip

\begin{defin} We say that a  continuous  function $v(x,t)$,
$v: \Gamma  \times \Rds^+\to\Rds$, is a {\em supersolution}   if
\begin{itemize}
   \item[{\bf (i)}] $v\circ \gamma$ is a viscosity supersolution of
   \eqref{HJg} in  $\Qcal$ for any arc $\gamma$;
     \item[{\bf (ii)}] for any vertex $x$ and  time $t_0 >0$,  if
     \[ \frac d{dt} \phi (t_0) < c_x\]
     for some $C^1$ subtangent $\phi$ to $v(x,\cdot)$ at $t_0$, then
there is  an arc $\gamma \in \Gamma_x$
such that  all the $C^1$  subtangents   $\Phi$,
constrained  to $\overline\Qcal$, to
$v \circ \gamma$  at  $(1,t_0)$   satisfy
     \[   \Phi_t(1,t_0) +  H_\gamma(1, \Phi'(1,t_0))  \geq 0.\]
\end{itemize}
\end{defin}

\smallskip

\begin{defin}
 We say that a  continuous  function $v(x,t)$, $v: \Gamma \times \Rds^+  \to \Rds$, is
 a {\em subsolution}  if
\begin{itemize}
   \item[{\bf (i)}] $v\circ \gamma$ is a viscosity subsolution  of \eqref{HJg} in
   $\Qcal$ for any arc $\gamma$;
     \item[{\bf (ii)}] for any vertex $x$ and  time $t_0 > 0$,
     all supertangents $\psi(t)$ to $v(x,\cdot)$ at $t_0$ satisfy
     \[  \frac d{dt} \psi(t_0) \leq c_x.\]
\end{itemize}
\end{defin}

We finally say that a  continuous function $v$ is {\em solution}
to (HJ$\Gamma$)   if it subsolution and supersolution at the same time.\\

We recall the following result of \cite{Siconolfi}:

\begin{theo}\label{orte}  For any continuous initial datum $u_0$ and flux limiter
$c_x$, there exists one and only one continuous solution
to (HJ$\Gamma$)  continuously attaining the datum $u_0$ at $t=0$.
\end{theo}
\smallskip

From now on we fix a continuous initial datum $u_0$ and a flux limiter $c_x$.

\medskip

\subsection{The equation on a single arc}

 We fix an arc $\gamma$,
  Given a continuous boundary datum
\[g: \partial\Qcal \to \Rds,\]
  we consider the equation \eqref{HJg}
coupled with the boundary conditions
\begin{equation}\label{boundary}
 U(s,t)= g(s,t) \qquad\hbox{for $(s,t) \in \partial\Qcal$.}
\end{equation}
or
\begin{equation}\label{boundary-}
 U(s,t)= g(s,t) \qquad\hbox{for $(s,t) \in \partial^-\Qcal$.}
\end{equation}
We say that the boundary data in \eqref{boundary} (resp. \eqref{boundary-}) are
{\em admissible} if there exists a continuous function  $V: \overline \Qcal \to \Rds$
solution of \eqref{HJg} in $\Qcal$ and equal to $g$ on $\partial \Qcal$
(resp. $\partial^-\Qcal$). The following results are well known, a proof is provided in
Appendix \ref{one}  for reader's convenience.

\begin{theo}\label{Q} If the datum $g$ in \eqref{boundary} is admissible then
the formula
\begin{equation}\label{Q1}
 V(s,t) = \inf \left \{ g(s_0,t_0) + \int_{t_0}^t L_\gamma(\eta, \dot\eta) \,d\tau \right \},
\end{equation}
where the infimum is over the elements $(s_0,t_0) \in \partial \Qcal$ with $t_0 <t$,
and the curves $\eta:[t_0,t] \to \overline \Qcal$ with $\eta(t_0)=s_0$,
 $\eta(t)=s$, is a solution to \eqref{HJg}, \eqref{boundary}.
\end{theo}

\smallskip

\begin{theo}\label{Q-} If the datum $g$ in \eqref{boundary-} is admissible then the  function
\begin{equation}\label{Q-1}
  W(s,t) = \inf \left \{ g(s_0,t_0) + \int_{t_0}^t L_\gamma(\eta, \dot\eta) \,d\tau \right\},
\end{equation}
where the infimum is over the elements $(s_0,t_0) \in \partial^- \Qcal$ with $t_0 <t$, and
the curves $\eta:[t_0,t] \to \overline \Qcal$ with $\eta(t_0)=s_0$,
 $\eta(t)=s$, is solution to \eqref{HJg}, \eqref{boundary-} and satisfies
\[\Psi_t(1,t)+ H_\gamma(1, \Psi'(1,t)) \geq 0\]
for any $t >0$, any subtangent $\Psi$ to $W$, constrained to $\overline \Qcal$, at $(1,t)$.
\end{theo}

\medskip

\subsection{The main result}

\smallskip
We  define for $(x,t) \in \Gamma \times \Rds^+$
\begin{equation}\label{repre}
u(x,t) = \inf_\xi \left \{\int_0^t L(\xi,\dot\xi) \, ds + u_0(\xi(0)) \right \},
\end{equation}
where $\xi$ is any curve from $[0,t]$ to $\Gamma$ with $\xi(t)=x$.

\smallskip

\begin{prop}\label{continuo} The function $u$ defined in \eqref{repre} is continuous
in $\Gamma \times \Rds^+$.
\end{prop}

\proof The continuity at any point $(x_0,t_0)$ with $ t_0 > 0$ comes straightforwardly from
the continuity of the minimal action functional and of $u_0$. Assume now that $(x_n,t_n)$
converges to $(x_0,0)$ and $\lim_n u(x_n,t_n) =:a$. We have
\[u(x_n,t_n) \leq u_0(x_n) + \int_0^{t_n} L(x_n,0) \, d\tau\]
which, sending $n$ to infinity, shows that $a \leq u_0(x_0)$. Assume now that
\begin{equation}\label{continuo1}
 u(x_n,t_n) = u_0(y_n)+ \int_0^{t_n} L(\xi_n,\dot\xi_n) \, d\tau,
\end{equation}
we claim that $y_n \to  x_0$. If not $y_n \to y_0\neq x_0$, up to
subsequences, and we find by the superlinearity of $L$  two sequences of real numbers
 $\alpha_k$, $\beta_k$ with $\beta_k \to + \infty$ such that
 \[L(x,q) \geq \alpha_k + \beta_k \, |q|.\]
 We derive that
 \[\int_0^{t_n} L(\xi_n,\dot\xi_n) \, d\tau \geq \alpha_k \, t_n + \beta_k \, d_\Gamma(x_n,y_n)
 \to + \infty \quad\hbox{as $k \to + \infty$, for any $n$}\]
 and
 \[\liminf_{n \to + \infty} \int_0^{t_n} L(\xi_n,\dot\xi_n) \, d\tau \geq \beta_k \,
 d_\Gamma(x_0,y_0) \quad\hbox{for any $k$.}\]
 This contradicts \eqref{continuo1}. Finally, we have
 \[u(x_n,t_n) \geq u_0(y_n) + t_n \, m_L,\]
 with $m_L$ defined as in \eqref{ml}, which implies $a \geq u_0(x_0)$,
 and concludes the proof.
\endproof

\smallskip

We keep in mind existence and uniqueness Theorem \ref{orte}. Our main result is:

\smallskip

\begin{theo}\label{main}  The function $u$ defined in \eqref{repre} is the
unique  solution to (HJ$\Gamma$)  with initial datum $u_0$ and flux
limiter $c_x$ for $x \in \Vbf$.
\end{theo}

We start by a preliminary result. We fix an arc $\gamma$.

\begin{prop}\label{propmain} Let $u$ be the function defined in \eqref{repre}, then
 the boundary datum
 \begin{equation}\label{propmain0}
  g(s,t) = \left \{ \begin{array}{cc}
             u_0\circ \gamma(s) & \quad\hbox{in $[0,1] \times \{0\}$} \\
             u\circ\gamma(0,t) & \quad\hbox{in $\{0\} \times \Rds^+$}\\
             u\circ\gamma(1,t) & \quad\hbox{in $\{1\} \times \Rds^+$},
           \end{array} \right .
 \end{equation}
where $u\circ\gamma$ is define in \eqref{eq:uconcg}, is admissible for \eqref{HJg} with boundary conditions \eqref{boundary}, and $u \circ\gamma$ is the
  solution given by  formula \eqref{Q1}.
\end{prop}
\proof We know that $u \circ \gamma$ is continuous in $\overline \Qcal$ by Proposition
\ref{continuo}. It is then enough by Theorem \ref{Q} to show that $u \circ \gamma$ is
given by \eqref{Q1}.
Let $(s,t) \in \Qcal$, we denote by
$\xi: [0,t] \to \Gamma$ an optimal curve for $u(\gamma(s),t)$, and set
\begin{eqnarray*}
  A &=&  \{\tau < t \mid \xi(\tau) \in \{\gamma(0), \gamma(1)\}\} \\
  t_0 & =&  \left \{\begin{array}{cc}
          \max A & \quad\hbox{if $A \neq \emptyset$} \\
           0 & \quad\hbox{if $A = \emptyset$}
         \end{array} \right .
\end{eqnarray*}
It is clear that $\xi((t_0,t]) \subset \gamma( (0,1))$. If $t_0 >0$, it is also
clear that the curve $\xi$ restricted to $[0,t_0]$ is optimal for
$u(\xi(t_0),t_0)$, hence we have by \eqref{lagrange}, Lemma \ref{deri}
and \eqref{deri1}
\begin{eqnarray*}
  u(\gamma(s),t) &=& u_0(\xi(0)) + \int_0^{t_0} L(\xi,\dot\xi) \, d\tau +
   \int_{t_0}^t L(\xi,\dot\xi) \, d\tau\\
   &=& u(\xi(t_0),t_0) +
   \int_{t_0}^t L_\gamma (\gamma^{-1}\circ\xi(\tau),d/d\tau (\gamma^{-1}\circ\xi(\tau)))
   \, d\tau\\
   &=& g(\gamma^{-1} \circ \xi(t_0),t_0) + \int_{t_0}^t L_\gamma (\gamma^{-1}\circ\xi(\tau),d/d\tau (\gamma^{-1}\circ\xi(\tau)))
   \, d\tau,
\end{eqnarray*}
 and the same formula holds true if $t_0=0$.
 Taking into account that $\gamma^{-1}\circ\xi$ is a curve
 by Lemma \ref{deri}, we derive that
 \[u \circ \gamma \geq V     \qquad\hbox{in $\Qcal$,}\]
 where $V$ is the function defined in \eqref{Q1}. Assume for purposes  of contradiction that
 there is $(s_1,t_1) \in \partial \Qcal$ with $t_1 <t$, and a curve
 $\eta: [t_1,t] \to \overline \Qcal$ with $\eta(t_1)=s_1$, $\eta(t)=s$ such that
 \begin{equation}\label{propmain1}
   u(\gamma(s),t) > g(s_1,t_1) + \int_{t_1}^t L_\gamma(\eta(\tau),\dot\eta(\tau)) \,d\tau.
 \end{equation}
 We denote by $\bar\xi:[0,t_1] \to \Gamma$ an optimal curve for $u(\gamma(s_1),t_1)$,
 we further
 denote by $\zeta:[0,t] \to \Gamma$ the curve obtained by concatenating $\bar\xi$ and
 $\gamma \circ \eta$, namely
 \[\zeta(\tau)= \left \{ \begin{array}{cc}
                  \bar\xi(\tau) & \quad \tau \in [0,t_1) \\
                 \gamma\circ \eta(\tau) & \quad \tau \in [t_1,t)
                \end{array} \right . \]
 We then deduce from \eqref{propmain1}
 \begin{eqnarray*}
   u(\gamma(s),t) &>& u(\gamma(s_1),t_1) +
   \int_{t_1}^t L_\gamma(\eta(\tau),\dot\eta(\tau)) \,d\tau \\
    &=& u_0(\bar\xi(0)) + \int_0^{t_1} L(\bar\xi,\dot{\bar\xi}) \, d\tau +
     \int_{t_1}^t L_\gamma(\eta(\tau),\dot\eta(\tau)) \,d\tau \\
     &=& u_0(\zeta(0)) + \int_0^t L(\zeta,\dot\zeta) \, d\tau,
 \end{eqnarray*}
 which cannot be. The proof is therefore concluded.
\endproof

\smallskip

\proof [{ Proof of Theorem \ref{main}}] We know, thanks to Proposition \ref{propmain}
that $u\circ \gamma$ is viscosity solution to \eqref{HJg} for any $\gamma$. It is then
 enough to check the solution conditions at the vertices. Assume by contradiction that
 there is $x \in \Vbf$,  $t_0>0$ and a supertangent $\psi$  to $u(x,\cdot)$ at $t_0$ with
 \[\frac d{dt} \psi(t_0) > c_x.\]
 We derive that  there exists $\delta >0$ with
 \[u(x,t_0-\delta) + h \, c_x \leq \psi(t_0-\delta) + \delta \, c_x < \psi(t_0)=u(x,t_0),\]
therefore
\[u(x,t_0) > u(x,t_0-\delta) + \int_{t_0-\delta}^{t_0} L(x,0) \, d\tau\]
which is in contrast with the very definition of $u$. We have so proved
the subsolution condition for $u$ at the vertices. We proceed considering an
admissible curve $\xi:[0,t_0] \to \Gamma$ realizing $u(x,t_0)$, we further  assume
the existence
of a subtangent  $\phi$ to $u(x,\cdot)$ at $t_0$ with
\[\frac d{dt} \phi(t_0) < c_x.\]
This implies that for any $\delta >0$ sufficiently small we have
\[u(x,t_0-\delta) + \delta \, c_x \geq \psi(t_0-\delta) + \delta \, c_x > \psi(t_0)=u(x,t_0).\]
We in turn derive that
\[\xi(t) \neq x   \qquad \hbox{for $\delta>0$ small, and $t \in (t_0-\delta,t_0)$.}\]
Taking into account the definition of admissible curve, there thus exists $0 \leq t_1 <t_0$,
$\gamma \in \Gamma_x$ such that $\xi((t_1,t_0)) \subset \gamma(\Qcal)$ and
\[u(x,t_0) =  \left \{\begin{array}{cc}
               u(\gamma(0),t_1) + \int_{t_1}^{t_0} L(\xi,\dot \xi) \, d\tau &
               \quad\hbox{if $t_1 >0$}  \\
               u_0(\gamma(s)) + \int_0^{t_0} L(\xi,\dot \xi) \, d\tau &
                \quad\hbox{for some $s \in [0,1)$ if $t_1 =0$}
             \end{array} \right .\]
If we define the boundary datum $g$ as in \eqref{propmain0} and set $\eta= \gamma^{-1} \circ \xi$,
 we can also write
\[u(\gamma(1),t) =  \left \{  \begin{array}{cc}
                      g(0,t_1) + \int_{t_1}^{t_0} L_\gamma(\eta,\dot\eta) \, d\tau&
                      \quad\hbox{if $t_1 >0$} \\
                      g(s,0) + \int_0^{t_0} L_\gamma(\eta,\dot\eta) \, d\tau &
                      \quad\hbox{for some $s \in [0,1)$ if $t_1 =0$}
                    \end{array} \right .\]
Since by Proposition \ref{propmain} $u \circ \gamma$ is given by formula \eqref{Q1}, with datum $g$ taken continuously
on $\partial \Qcal$, we therefore see that
\[u \circ \gamma = W \quad\hbox{at $(1,t_0)$}, \quad u \circ \gamma \leq W
\quad\hbox{in $\overline Q$,}\]
where $W$ is defined as in Theorem \ref{Q-}. If $\varphi$ is subtangent, constrained to $\overline Q$,
 to $u\circ\gamma$  at $(1,t_0)$, then $\varphi$ is also
 subtangent to $W$ at $(1,t_0)$ and by Theorem \ref{Q-} we get
 \[\varphi_t(1,t_0) + H_\gamma(1,\varphi'(1,t_0)) \geq 0.\]
 The proof is then complete.
\endproof

\begin{appendix}

\section{Proof of Theorems \ref{Q}, \ref{Q-}} \label{one}

The proof is based on the following lemmata:

\begin{lem}\label{subsol}
Let $U:\overline\Qcal\to\Rds$ be a continuous function.
If for any $(s,t)\in\Qcal$, one has
\[U(s,t)-U(\eta(t-\delta), t-\delta)\le\int_{t-\delta}^t L_\gamma(\eta,\dot\eta) \, d\tau\]
for any $C^1$ curve $\eta$ taking values
in $[0,1]$ with $\eta(t)=s$, $\delta >0$ small enough,
 then $U$ is a subsolution to \eqref{HJg}.
\end{lem}
\proof
We fix $\lambda \in \Rds$. Let  $\Phi$ be a supertangent to $U$ at $(s,t)$.
If $\eta$ is a
 $C^1$ curve with $\eta(t)=s$ and $\dot\eta(t)=\lambda$, we have for any $\delta>0$ small enough
\[
\frac{\Phi(s,t)-\Phi(\eta(t-\delta),t-\delta)}\delta\le\frac{U(s,t)-U(\eta(t-\delta),t-\delta)}
\delta\le\frac1\delta\int_{t-\delta}^tL_\gamma(\eta,\dot\eta)\, d\tau.
\]
Sending $\delta$ to 0 we get
\[\Phi_t(s,t)+  \lambda \, \Phi'(s,t) - L_\gamma(s,\lambda)\le0.
\]
By the very definition of $L_\gamma$ and the fact that
$\lambda$ has been taken  arbitrarily, we further get
\[\Phi_t(s,t)+H_\gamma(s,\Phi'(s,t))\le 0,\]
which gives the assertion.
\endproof

\smallskip

\begin{lem}\label{supsol}
Let $U:\overline\Qcal\to\Rds$ be a continuous function.
If for any $(s,t)\in\Qcal$ there is a curve $\eta$ taking values
in $[0,1]$ with $\eta(t)=s$  such that
\[U(s,t)-U(\eta(t-\delta), t-\delta)\ge\int_{t-\delta}^t L_\gamma(\eta,\dot\eta) \, d\tau,
 \]
when  $\delta >0$ is small enough, then $U$ is a supersolution to \eqref{HJg} in $\Qcal$.
\end{lem}
\proof
Let  $\Phi$ be a subtangent to $U$ at $(s,t)$,
then we have  for any  $\mu \in \Rds$, $\delta$ small
\begin{eqnarray}
\frac{\Phi(s,t)-\Phi(\eta(t-\delta), t-\delta)}\delta &\geq&\frac{U(s,t)-
U(\eta(t-\delta), t-\delta)}\delta \nonumber\geq \frac1\delta\int_{t-\delta}^tL_\gamma(\eta,\dot\eta) \, d\tau\\
&\geq&\frac1\delta\int_{t-\delta}^t (\mu \, \dot\eta(\tau)-
H_\gamma(\eta(\tau),\mu)) \, d\tau  \label{eq:supsol1}\\
&=&\, \mu  \,  \frac{\eta(t)-\eta(t-\delta)}\delta-\frac1\delta\int_{t-\delta}^t
H_\gamma(\eta(\tau),\mu) \,d\tau. \nonumber
\end{eqnarray}
Given any infinitesimal sequence $\delta_n$,  we find,
by the mean value theorem applied to the function $\Phi$, points $(s_n,t_n)$ in the
segment joining
$(\eta(t-\delta_n),t-\delta_n)$ to $(s,t)$ with
\[
\Phi(s,t)-\Phi(\eta(t-\delta_n),t-\delta_n)=
\Phi_t(s_n,t_n)\, \delta_n+\Phi'(s_n,t_n)(\eta(t)-\eta(t-\delta_n)),
\]
therefore \eqref{eq:supsol1} with $\mu = \Phi'(t_n,s_n)$ yields, for any $n\in\Nds$,
\begin{multline*}
\Phi_t(t_n,s_n)+ \Phi'(t_n,s_n)\frac{\eta(t)-\eta(t-\delta_n)}{\delta_n}\\
\ge \Phi'(t_n,s_n)\frac{\eta(t)-\eta(t-\delta_n)}{\delta_n}-\frac1{\delta_n}
\int_{t-\delta_n}^t H_\gamma(\eta(\tau),\Phi'(s_n,t_n)) \, d\tau.
\end{multline*}
We thus have
\[
\Phi_t(s,t)+H_\gamma(s,\Phi'(s,t))=\lim_{n\to\infty}\Phi_t(t_n,s_n)
+\frac1{\delta_n}\int_{t-\delta_n}^t H_\gamma(\eta,\Phi'(t_n,s_n)) \, d\tau\ge0,
\]
 This shows the assertion, since $(s,t)$ and $\Phi$ are arbitrary.
\endproof

\smallskip

\begin{lem}\label{locsolcont}
The value functions defined in \eqref{Q1} and \eqref{Q-1} are continuous.
\end{lem}
\proof
We will only prove that \eqref{Q1} is continuous, since the proof for \eqref{Q-1}  is identical.\\
We consider  the set
\[ E = \{ (s_0,t_0,s,t)\in\overline\Qcal^2 \mid t >t_0 \;\hbox{or} \; t=t_0, \, s=s_0\}\]
and  the function
\[S_\gamma(s_0,t_0,s,t):=\inf\left\{\int_{t_0}^tL_\gamma(\eta,\dot\eta)d\tau \mid \eta\text{ curve with }\eta(t_0)=s_0,\eta(t)=s\right\}\]
in  $E$, and $S_\gamma=+\infty$ in $\overline \Qcal \setminus E$.
It is continuous  when $t > t_0$ and lower semicontinuous  otherwise, see \cite{fathi14} and the proof of Proposition \ref{continuo}. This implies that if $\{(s_n,t_n)\}_{n\in\Nds}$ is a  sequence in $\overline\Qcal$ converging to an element  $(\bar s,\bar t)$,  there is  $(s_0,t_0)\in\partial\Qcal$ with $t_0 \leq \bar t$ and a sequence $\{(s_n',t_n')\}_{n\in\Nds}\subset\partial\Qcal$ with $t'_n \leq t_n$  such that
\[
V(\bar s,\bar t)=g(s_0,t_0)+S_\gamma(s_0,t_0,\bar s,\bar t),\qquad V(s_n,t_n)=g(s_n',t_n')+S_\gamma(s_n',t_n',s_n,t_n).
\]
When $\bar t>t_0$ the continuity is easily checked,  we  therefore  focus on the $\bar t=t_0$, $\bar s=s_0$ case. In other term, we assume $(\bar s,\bar t) \in \partial\Qcal$ and, taking into account that $g$ is admissible,  $V(\bar s,t)= g(\bar s,\bar t)$.  We further assume that $V(s_n,t_n)$ converges to a quantity $a$. The lower semicontinuity of $S_\gamma$ and the admissibility of $g$ yield
 \begin{equation}\label{locsolcont1}
 a=\lim_n g(s_n',t_n')+S_\gamma(s_n',t_n',s_n,t_n)
\geq g(\bar s, \bar t)= V(\bar s, \bar t).
 \end{equation}
 If $t_n =t'_n$, up to a subsequence, then $(s_n,t_n) \in \partial\Qcal$ and  $V(s_n,t_n)=g(s_n,t_n)$ so that
 \begin{equation}\label{locsolcont2}
 a= \lim_n g(s_n,t_n)= g(\bar s,\bar t)=V(\bar s,\bar t).
 \end{equation}
If $\bar t=0$ then
\begin{equation}\label{locsolcont3}
  a \leq  \lim_n g(s_n,0) +t_n \, \max_{s\in [0,1]} L_\gamma(s,0) = g(\bar s,0)=V(\bar s,\bar t).
\end{equation}
If $\bar t>0$  and  $t'_n< t_n$, for any $n$, then we can choose $t''_n  >0$  such that
\[t''_n < t_n, \quad |s_n-\bar s| \leq t_n - t''_n, \quad t''_n \to \bar t\]
so that we have
\begin{equation}\label{locsolcont4}
 a \leq  \lim_n g(\bar s,t''_n) + \max_{(s,\lambda)\in[0,1]\times[-1,1]}L_\gamma(s,\lambda) \, (t_n-t''_n) = g(\bar s,\bar t)=V(\bar s, \bar t).
\end{equation}
By combining \eqref{locsolcont2}, \eqref{locsolcont3}, \eqref{locsolcont4} with \eqref{locsolcont1}, we get the assertion.
\endproof

\smallskip

\proof [Proof of Theorem \ref{Q}] We take for granted that for any $(s,t) \in \Qcal$
there exists  a curve $\eta$ realizing the infimum in formula \eqref{Q1}. Then we have for
 any $\delta >0$ small enough
 \[V(s,t)= V(\eta(t-\delta),t-\delta) + \int_{t-\delta}^t L_\gamma(\eta,\dot\eta) \, d\tau\]
 which implies by Lemma \ref{supsol} that $V$ is a subsolution. Since $V$ is defined by
 \eqref{Q1} and is continuous by Lemma \ref{locsolcont}, it also satisfies the assumption of Lemma \ref{subsol}, which shows that it
 is subsolution  as well.
\endproof

\smallskip

\proof [Proof of Theorem \ref{Q-}] Arguing as in the proof of Theorem \ref{Q}, we see that
the function defined in \eqref{Q-1} is solution to \eqref{HJg}. If we take a point $(1,t)$ for
some $t>0$, then there exists an optimal curve, say $\eta$, realizing $W(1,t)$,
 if we consider
a subtangent $\Phi$, constrained to $\overline \Qcal$,  to $W$ at $(1,t)$ then arguing as in
Lemma \ref{supsol} we derive that
\[\Phi_t(1,t) + H_\gamma(1,\Phi'(1,t)) \geq 0.\]
This concludes the proof.
\endproof

\bigskip

\section{Lipschitz regularity of the minimizers} \label{covidcovid}
We  present the results of this appendix in a rather  general frame.  They are applied in Theorem \ref{lippo}   to the one--dimensional Hamiltonians/Lagrangians appearing  in the equations on networks.

We consider an Hamiltonian $\wtd H: \Rds^N \times \Rds^N \to \Rds$ enjoying the usual properties of continuity in both arguments plus convexity and superlinearity in the second one, and denote by $\wtd L: \Rds^N \times \Rds^N \to \Rds$ the corresponding Lagrangian defined through Fenchel transform.   We provide  the notion of $a$--Lagrangian parametrization of a curve.

\smallskip

\begin{defin} Given  an absolutely continuous  curve $\eta:[0,T] \to \mathbb R^N$, we say that $a \in \Rds$ is {\em admissible} if
\[a \geq \max_{x \in \spt \eta}\,\min_{p \in\Rds^N} \wtd H(x,p).\]
 We say that $\eta$ has an
 {\em $a$--Lagrangian parametrization}, for  some admissible $a$,  if
\[\wtd L(\eta(t),\dot\eta(t)) = \sigma_a(\eta(t),\dot\eta(t)) -a  \qquad\hbox{for a.e. $t \in [0,T]$,}\]
where
\[\sigma_a(x,q)= \max \{ p \cdot q \mid p \;\hbox{with $\wtd H(x,p) \leq a$}\}.\]
Note that  the sublevels in the above formula are nonempty because of the admissibility of $a$.
\end{defin}

\smallskip

\begin{prop} Any curve  $\eta: [0,T] \to \mathbb R^N$  with an $a$--Lagrangian parametrization is Lipschitz continuous.
\end{prop}
\begin{proof}   By assumption we have
\[\wtd L(\eta(t), \dot\eta(t)) = p(t) \, \dot\eta(t) - \wtd H(\eta(t),p(t)) \qquad\hbox{for a.e. $t$,}\]
where $p(t)$ satisfies $\wtd H(\eta(t),p(t))=a$. This implies that
$\dot\eta(t) \in \partial_p \wtd H(\eta(t),p(t))$, and  by the uniform coercivity of $\wtd H$ in $p$, we get that $|p(t)| < M$, for a.e. $t$ and some $M >0$. Since $\wtd H$ is locally Lipschitz continuous in $p$ uniformly in $\spt \eta$, see \eqref{newnew}, we find a constant $C_M$ with $|\dot\eta(t)|< C_M$ for a.e. t, concluding the proof.
\end{proof}

\medskip

The main result  we aim at showing is:

\begin{theo}\label{covid} Given an absolutely continuous curve $\xi: [0,T] \to \Rds^N$, the problem
\[ \inf \left \{\int_0^T \wtd L(\zeta,\dot\zeta) \, dt \mid \zeta:[0,T] \to \Rds^N \;\hbox{A.C. with}\; \xi(0)=\zeta(0),\,\xi(T)= \zeta(T),\, \spt \zeta \subset \spt \xi  \right\}  \]
admits minimizers. All of them possess  an $a$--Lagrangian representation, for some admissible $a$,  and
consequently are Lipschitz continuous.
\end{theo}

\smallskip

We need  a further notion.

\begin{defin} Given an absolutely continuous curve $\xi:[0,T] \to \Rds^N$, a curve $\zeta:[0,T'] \to \Rds^N$ is said a {\em reparametrization} of $\xi$ if there exists a nondecreasing absolutely continuous function $\varphi$ from $[0,T']$ onto $[0,T]$ with
\[ \zeta(t)= \xi \circ \varphi(t) \qquad\hbox{for any $t \in [0,T']$.}\]
\end{defin}

\smallskip
Since the composition of two absolutely continuous functions, with the second one monotonic, is still absolutely continuous, we see that any reparametrization $\zeta$ of $\xi$ keeps absolutely continuity. It is further immediate that the support  of a curve and its length are invariant under reparametrization.  Similarly,  exploiting the positive homogeneity of $\sigma_a(x,\cdot)$ we see that
\begin{equation}\label{precovid}
 \int_0^{T'} \sigma_a(\zeta,\dot\zeta) \, dt = \int_0^T \sigma_a(\xi,\dot\xi) \, dt
\end{equation}
for any admissible $a$. Note finally  that   if $\zeta$ is  a reparametrization of $\xi$, the converse property  in general  is not true  for $\varphi$ could not have strictly positive derivative for a.e. $t$, see Zarecki criterion for an absolutely continuous inverse in \cite{Bernal}. We have:

\begin{prop}\label{finecovid} Any absolutely continuous curve  in $[0,T]$  is a reparametrization of a curve $\eta$ with constant speed, namely absolutely continuous with $|\dot\eta(t)| \equiv \mathrm{constant}$ for a.e. $t$, defined in the same interval.
\end{prop}
\proof
See \cite{davini07,Bernal}.
\endproof

\smallskip

The following crucial result of \cite{davini07} is in the spirit of Erdmann condition.

\begin{theo}\label{davini} Given a curve $\xi:[0,T] \to \Rds^N$  with constant speed, there is a reparametrization in the same interval with an $a$--Lagrangian parametrization, for some admissible $a$.
\end{theo}
\begin{proof}
See Theorem 3.16  and Remark 3.17 in \cite{davini07}.
\end{proof}

\begin{proof}[{\it Proof of Theorem \ref{covid}}] \,
By the usual existence argument, see Theorem \ref{funda}, and the fact that $\spt \xi$ is compact,  we find a minimizer  $\zeta_0:[0,T] \to \Rds^N$  of the optimization problem in the statement. By Proposition \ref{finecovid}, $\zeta_0$ is the reparametrization of curve $\eta_0:[0,T] \to \Rds^N$ with constant speed. Taking into account \eqref{precovid}, we have
\begin{equation}\label{covid1}
 \int_0^T \wtd L(\zeta_0,\dot\zeta_0) \, dt \geq \int_0^T \sigma_a(\zeta_0,\dot\zeta_0) \, dt - a \, T = \int_0^T \sigma_a(\eta_0,\dot\eta_0) \, dt - a \, T
\end{equation}
for any admissible $a$. By Theorem \ref{davini}, there is an admissible $a$  and a reparametrization $\zeta$ of $\eta_0$ with an $a$--Lagrangian parametrization. We therefore get
\[\int_0^T \wtd L(\zeta,\dot\zeta) \, dt = \int_0^T \sigma_a(\zeta,\dot\zeta) \, dt - a \, T = \int_0^T \sigma_a(\eta_0,\dot\eta_0) \, dt - a \, T \leq \int_0^T \wtd L(\zeta_0,\dot\zeta_0) \, dt.\]
We deduce from the optimality of $\zeta_0$ that equality must prevail in the above formula, as well as in \eqref{covid1}. This implies the assertion.
\end{proof}

\end{appendix}

\printbibliography[heading=bibintoc]
\end{document}